\documentclass[11 pt]{article}
 \usepackage{amssymb}
\usepackage{verbatim}
\usepackage{latexsym}
\usepackage{eucal}
\newtheorem{theorem}{Theorem}[section]

\newtheorem{lemma}{Lemma}[section]

\newtheorem{remark}{Remark}[section]

\def\proof{\mbox {\bf Proof.}\quad}
\def\vfi{\varphi}
\newtheorem{sub}{\name}[section]
\newcommand{\bs}{
\begin{sub}}
\newcommand{\es}{
\end{sub}}
\newcommand{\bsl}[1]{
\begin{sub}\label{#1}}
\newcommand{\bth}[1]{\def\name{Theorem}
\begin{sub}\label{t:#1}}
\newcommand{\blemma}[1]{\def\name{Lemma}
\begin{sub}\label{l:#1}}
\newcommand{\bcor}[1]{\def\name{Corollary}
\begin{sub}\label{c:#1}}
\newcommand{\bdef}[1]{\def\name{Definition}
\begin{sub}\label{d:#1}}
\newcommand{\bprop}[1]{\def\name{Proposition}
\begin{sub}\label{p:#1}}
\newcommand{\brem}[1]{\def\name{Remark}
\begin{sub}\label{r:#1}}


\newcommand{\BA}{
\begin{array}}
\newcommand{\EA}{
\end{array}}
\newcommand{\BAN}{\renewcommand{\arraystretch}{1.2}
\setlength{\arraycolsep}{2pt}
\begin{array}}
\newcommand{\BAV}[2]{\renewcommand{\arraystretch}{#1}
\setlength{\arraycolsep}{#2}
\begin{array}}
\newcommand{\BSA}{
\begin{subarray}}
\newcommand{\ESA}{
\end{subarray}}
\newcommand{\BAL}{
\begin{aligned}}
\newcommand{\EAL}{
\end{aligned}}
\newcommand{\BALG}{
\begin{alignat}}
\newcommand{\EALG}{
\end{alignat}}
\newcommand{\BALGN}{
\begin{alignat*}}
\newcommand{\EALGN}{
\end{alignat*}}

\newcommand{\qed}{\\
${}$ \hfill $\square$}

\newcommand{\forevery}{\quad \forall}

\newcommand{\abs}[1]{\left |#1\right |}

\def\angb<#1>{\langle #1 \rangle}

\newcommand{\myfrac}[2]{{\displaystyle \frac{#1}{#2} }}
\newcommand{\myint}[2]{{\displaystyle \int_{#1}^{#2}}}

\newcommand{\prt}{
\partial}

\newcommand{\ti}{\times}

\newcommand{\nind}{\noindent}

\def\ga{\alpha}            \def\gg{\gamma}
       \def\gd{\delta}      \def\ge{\epsilon}

\def\gm{\mu}                 
    \def\gr{\rho}        
\def\gs{\sigma}       \def\gt{\tau}

\def\Gg{\Gamma}     \def\Gd{\Delta}

   \def\CB{{\mathcal B}}

   \def\BBR {\mathbb R}

\catcode`\@=11

\def\la{\lambda}
\def\vep{\varepsilon}
\def\pd#1#2{\frac{\partial #1}{\partial #2}}
\def\eqalign#1{\null\,\vcenter{\openup1\jot \m@th
  \ialign{\strut\hfil$\displaystyle{##}$&$\displaystyle{{}##}$\hfil
     &&\strut$\displaystyle{##}$\hfil&$\displaystyle{{}##}$
     \hfil\crcr#1\crcr}}\,}
\newcommand{\rife}[1]{(\ref{#1})}
\def\qed{{\unskip\nobreak\hfil\penalty50
          \hskip2em\hbox{}\nobreak\hfil\mbox{\rule{1ex}{1ex} \qquad}
   \parfillskip=0pt
   \finalhyphendemerits=0\par\medskip}}

\begin{document}
\title {$\qquad$\bf Symmetry of large solutions of nonlinear elliptic 
equations in a ball}
\author{{\bf\large Alessio Porretta}
\footnote{The author acknowledges the support of RTN european project:
FRONTS-SINGULARITIES, RTN contract: HPRN-CT-2002-00274. }
\hspace{2mm}\vspace{2mm}\\
{\it Dipartimento di Matematica,  Universit\`a di Roma Tor 
Vergata},\\
{\it Via della Ricerca Scientifica 1, 00133 Roma, Italia}\\
\vspace{2mm}\\
{\bf\large Laurent V\'eron}
\vspace{2mm}\\
{\it Laboratoire de Math\'ematiques et Physique Th\'eorique
 CNRS UMR  6083},
\\
{\it Universit\'e Fran\c{c}ois Rabelais},{\it Tours 37200, France}}

\date{}
\maketitle

{\small{\bf Abstract} Let $g$ be a locally Lipschitz continuous real 
valued function which satisfies the Keller-Osserman condition and is 
convex at infinity, then any large solution of $-\Delta u+g(u)=0$ in 
a ball is radially symmetric.}
\vspace{1mm}
\hspace{.05in}

\section {Introduction}
\setcounter{equation}{0} 
Let $B_{R}$ denote the open ball of center $0$ and radius $R>0$ in $\BBR^N$, $N\geq 
2$. A classical result due to Gidas,  Ni and Nirenberg  \cite {GNN} 
asserts that, if $g$ is a locally Lipschitz continuous real valued 
function, any $u\in C^2(\overline \Omega)$ which is  a positive solution of 
\begin {equation}\label {simple}\left\{\BA {l}
-\Gd u+g(u)=0
\quad \mbox {in } B_{R}\\
\phantom{-\Gd u+gu}
u=0\quad \mbox {on } \prt B_{R}
\EA\right.
\end {equation} 
is radially symmetric. The proof of this result is based on the 
celebrated Alexandrov-Serrin moving plane method. Later on, this 
method was used in many occasions, with a lot of refinements for 
obtaining selected symmetry results and a priori estimates for 
solutions of semilinear elliptic equations. If the boundary condition 
is replaced by $u=k\in\BBR$, clearly the radial symmetry still holds 
if $u-k$ does not change sign in $B_{R}$. Starting from this 
observation, it was conjectured by Brezis \cite {Br} that any  solution $u$ of 
\begin {equation}\label {large}\left\{\BA {l}
-\Gd u+g(u)=0\quad \mbox {in } B_{R}\\
\phantom{-}\displaystyle \lim_{\abs{x}\to R }u(x)=\infty,
\EA\right.
\end {equation} 
is indeed radially symmetric. Notice that this problem admits a 
solution (usually called a  \lq\lq large solution\rq\rq) if and only if $g$ satisfies the Keller-Osserman condition: 
$g\geq h$ on $[a,\infty)$, for some  $a>0$ where $h$ is non 
decreasing and satisfies 
\begin {equation}\label {KO}\left\{\BA {l}
\myint {a}{\infty}\myfrac {ds}{\sqrt {H(s)}}<\infty\\
\mbox {where }H(s)=\myint{a}{s}h(t)dt.
\EA\right.
\end {equation} 
Up to now, at  least to our knowledge,  only partial results were known concerning  the radial symmetry of solutions of \rife{large}: in \cite{MRW}, the authors prove this result assuming (besides the Keller--Osserman condition) that $g'(s)/\sqrt{G(s)}\to \infty$ as $s\to \infty$, or for the special case when $g(s)=s^q$, using 
the estimates for the second term of the asymptotic expansion of the solution near the boundary. Of course, 
the symmetry can also be obtained via  uniqueness, 
however uniqueness is known under an assumption of {\emph {global}} monotonicity and 
convexity (\cite {MV}, \cite {MV1}). Otherwise, it is easy to prove,  
by a one--dimensional topological argument, that uniqueness for problem \rife{large} holds for almost all 
$R>0$ under the mere  monotonicity assumption. 
However, if $g$ is not 
monotone, uniqueness may not hold (see e.g.   \cite {ADL}, \cite{MRW}, \cite {Po}), and it turns out to be very important to know whether {\it all} the solutions constructed in  a ball are radially symmetric, a  fact that would lead to a  full classification of all possible solutions. Let us point out that the interest in such qualitative properties of large solutions has being raised in the last few years from different problems (see e.g. \cite {ADL}, \cite{DG}, \cite{DG1}, \cite{DY} and the references therein).

In this article we prove that  Brezis' conjecture is verified under an assumption of 
asymptotic convexity upon $g$, namely we prove

 \begin{theorem}\label{T}{\it  Let $g$ be a  locally Lipschitz continuous function. Assume that $g$ is positive and convex on 
$[a,\infty)$  for some $a>0$, and satisfies 
the Keller-Osserman condition. Then any $C^2$ solution of (\ref {large})  is 
radially symmetric and increasing}.
\end{theorem}

Notice that  the Keller-Osserman condition implies that the function 
$g$ is superlinear at infinity. The convexity assumption on $g$ is 
then very natural in such context. 

In order to prove Theorem \ref{T},  we prove  first  a suitable adaptation of Gidas-Ni-Nirenberg moving-planes method 
to the framework of large solutions, without requiring any monotonicity 
assumption on $g$. This first result, which can have an interest in its own, reads as follows:
 
 \begin{theorem}  {\it Assume that $g$ is locally Lipschitz continuous  and let $u$ be a solution of (\ref{large}) 
which satisfies
 \begin {equation}\label {cond}\left\{\BA {l}
\displaystyle \lim_{\abs{x}\to R}\prt_r u=\infty\\[4mm]
|\nabla_{\gt}u|=\circ\left(\prt_r u\right) 
\quad\mbox {as }\abs{x}\to R, \forevery \gt\perp x\;\mbox { s. t. } 
\abs \gt=1,
\EA\right.
\end {equation}
where $\prt_r u$ and $\nabla_{\gt}u$ are respectively the radial derivative
and the tangential gradient of $u$. Then $u$ is radially symmetric and 
$\prt_r u>0$ on 
$B_{R}\setminus\{0\}$.} 
\end{theorem}

Thus, in view of the previous statement, our main point in order to deduce the general result of Theorem \ref{T}  
is to  prove  that condition (\ref {cond}) always holds (even in a  stronger form) if  we assume that $g$ is asymptotically convex, and this is achieved by 
providing sharp informations on the radial and the tangential behavior of $u$ 
near the boundary.

\section{Proof of the results}
\setcounter{equation}{0} 

Let $\CB=\{e_{1},...,e_{N}\}$ be the canonical basis of $\BBR^N$. If 
$P\in\BBR^N$ and $\gr>0$, we denote by 
 $B_\gr(P)$ the open ball with center $P$ and radius $\gr$, and for 
simplicity $B_\gr(0)=B_{\gr}$. We consider the problem
\begin {equation}\label {large2}\left\{\BA {l}
-\Gd u+g(u)=0\quad \mbox {in } B_{R}\\
\phantom{----}\!u(x)=\infty\quad \mbox {on } \prt B_{R},
\EA\right.
\end {equation}
where $R>0$. By a solution of (\ref {large2}), we mean that $u\in 
C^2(B_R)$ is a classical solution in the interior of the ball and that
$u(x)$ tends to infinity uniformly as $|x|$ tends to $R$.
\vskip1em
We shall consider the following assumptions on $g$:
\begin{equation}\label{g1}
g\,:\,\BBR \to \BBR \qquad \hbox{is locally Lipschitz continuous}. 
\end{equation}
\begin{equation}\label{g2}
\exists\,\, a >0\,\,\,\hbox  {s.t.}\quad g \,\, \hbox{ is positive and  convex 
on } [a,\infty), 
\end{equation}
 and satisfies
\begin{equation}\label{g3}
\int_{a}^{+\infty} \myfrac{1}{\sqrt{G(t)}} dt 
<+\infty\,,\qquad\mbox{where }G(t)=\int_a^t g(s)ds.
\end{equation} 
 \medskip

Note that convexity and \rife{g3} imply that $g$ is increasing on 
$[b,\infty)$ for some $b>0$.\\

If $u\in C^1(B_{R})$ we denote by $\prt u/\prt r (x) =\langle Du(x), 
x/\abs x\rangle$ the radial derivative of $u$, and by 
$\nabla_{\gt}u(x)=(Du(x)-|x|^{-1}\,\prt u/\prt r (x))x$ the tangential 
gradient of $u$.
Our first technical result, which is a reformulation in the framework 
of large solutions of the famous original proof of Gidas, Ni and 
Nirenberg \cite{GNN}, is the following

\begin{theorem}\label{GNN}
Assume that $g$ satisfies (\ref {g1}), and let $u$ be a solution of 
(\ref {large2}). If there holds
\begin{equation}\label{der}\BA {l}
(i)\quad \displaystyle \lim_{|x|\to R}\myfrac {\prt u}{\prt r} (x) 
=\infty \\
(ii)\quad\abs {\nabla_{\gt}(x)}= \circ\left(\myfrac {\prt u}{\prt r} 
(x)\right)\quad\hbox{as } |x|\to R, 
 \EA
\end{equation}
then $u$ is radially symmetric and $\prt u/\prt r >0$ in $B_{R}\setminus \{0\}$.
\end{theorem}
\proof  Since the equation is invariant by rotation, it is 
sufficient to prove that (\ref{der})  implies that  $u$ is symmetric 
in the $x_1$ direction.\smallskip

\nind We claim first that for any $P\in \partial B^+\,:=\partial B_R\cap 
\{x\in \BBR^N\,:\,x_1>0\}$, there exists $\delta \in (0,R)$ such that
\begin{equation}\label{claim}
 \myfrac{\partial u}{\partial x_1}(x) >0\qquad \forall\,x\in B_R \cap 
B_{\gd}(P).
\end{equation}
Indeed, thanks to (\ref{der}) we have,
$$
\eqalign{ 
\myfrac{\partial u}{\partial x_1}&= \myfrac{\partial u}{\partial r} 
\myfrac{x_1}{|x|} + (D u- \myfrac{\partial u}{\partial r} 
\myfrac{x}{|x|})\cdot e_1 \cr &=  \myfrac{\partial u}{\partial 
r}\left( \myfrac{x_1}{|x|}+ \left( \myfrac{\partial u}{\partial 
r}\right)^{-1}\left(D u- \myfrac{\partial u}{\partial r} 
\myfrac{x}{|x|}\right)\cdot e_1\right)\cr & = \myfrac{\partial 
u}{\partial r}\left( \myfrac{x_1}{|x|}+ \circ(1)\right) \quad\hbox{as 
$|x|\to R$.}\cr}
$$
 Since $P\in \partial B^+$, the claim follows straightforwardly.
\medskip

Next we follow the construction in \cite{GNN}. For any $\la<R$, set 
$T_\la$ the hyperplane $\{x_1=\la\}$ and $\Sigma_\la=\{x\in B_R\,:\, 
\la<x_1<R\}$, $\Sigma'_\la= \{x\in B_R
\,:\, 2\la-R<x_1<\la\}$ the symmetric 
caps reflected with respect to $T_\la$; denote also 
$x_\la=(2\la-x_1,x_2,\ldots,x_N)$ the reflected point and 
$u_\la=u(x_\la)$ the reflected function, for $x\in \Sigma_\la$.  Let 
$P_{0}=Re_{1}$ and let $\delta= \delta (P_0)>0$ be the real number 
such that (\ref {claim}) holds in  $B_R\cap B_{\gd_{0}}(P_{0})$.
If $\la_0=R-\gd_{0}^2/2R$, there holds
$$
\myfrac{\partial u}{\partial x_1}>0 \quad \hbox {in}\quad 
\Sigma_{\la_0}\cup \Sigma'_{\la_0},
$$ 
so that, in particular, 
\begin{equation}\label{lam}
 u(x_\la)<u(x)\qquad \mbox{ and }\qquad \myfrac{\partial u}{\partial 
x_1} >0\quad\hbox{ in $\Sigma_\la$},
\end{equation}
for $\la\geq\la_0$. We define 
$$\mu=\inf\{\la>0\,:\, \mbox { s. t. }\rife{lam} \mbox { holds true } 
\}$$
 and we claim that $\mu=0$.
 We proceed by contradiction and assume that $\mu>0$. Denote by 
$K_\mu=T_\mu\cap \partial B_R$: since $K_\mu$ is compact, thanks to 
\rife{claim} there exists an $\vep$-neighborhood  $U_\vep$ of $K_\mu$ 
such that  
\begin{equation}\label{ue}
\myfrac{\partial u}{\partial x_1} >0 \qquad \mbox{in $U_\vep \cap 
B_R$}. 
\end{equation}
By definition of $\mu$ there holds $u\geq u_\mu$ in $\Sigma_\mu$; 
thus, if we denote $D_\vep=B_{R-\ge/2}\cap \Sigma_\mu$, 
we have
\begin{equation}\label{eps}
\cases{\Delta (u-u_\mu)= a(x) (u-u_\mu) &in $D_\vep$\cr
\phantom {--}u-u_\mu \geq 0 & in $D_\vep$,\cr}
\end{equation}
where $a(x)= (g(u)-g(u_\mu)/(u-u_\mu)$. Thanks to  \rife{g1} and 
since $\overline D_\vep$ is in the interior of $B_R$, $a(x)$ is a 
bounded function in $D_\vep$, and  the strong maximum principle 
applies to \rife{eps}. Since $u$ tends to infinity at the 
boundary and is finite in the interior, for $\vep$ small we clearly have  $u\not \equiv u_\mu$ in $D_\vep$:
therefore we conclude that  $u>u_\mu$ in $D_\vep$, and, since $u=u_\mu$ on $T_\mu\cap  \prt D_\vep$ and 
$\partial u_\mu/\partial 
x_1=-\partial u/\partial x_1$ on $T_{\gm}$, it follows from Hopf 
boundary lemma that
$$
\myfrac{\partial u}{\partial x_1} >0 \quad \hbox{ on $T_\mu\cap 
\partial  D_\vep$.}
$$
Since $u\in C^1(B_R)$, the last assertion, together with \rife{ue}, 
implies that there exists $\sigma>0$ such that
\begin{equation}\label{sig}
\myfrac{\partial u}{\partial x_1} >0 \quad \hbox{in $B_R\cap \{x\,:\, 
\mu-\sigma<x_1< \mu+\sigma\}$.}
\end{equation}
Moreover, since $\vep$ can be chosen arbitrarily  small, we deduce 
that 
$$
u>u_\mu \quad \hbox{in $\Sigma_\mu$.}
$$
Now, by definition of $\mu$,  there exists  an increasing sequence 
$\la_n$ converging to $\mu$ and points $x_n\in \Sigma_{\la_n}$ such 
that 
\begin{equation}\label{cont}
u(x_n)\leq u((x_n)_{\la_n}).
\end{equation}
Up to subsequences, $\{x_n\}$ will converge to a point $\bar x\in 
\overline \Sigma_\mu$. However, $\bar x$ cannot belong to 
$\Sigma_\mu$, since in the limit  we would have $u(\bar x)\leq u(\bar 
x_\mu)$ while we  proved that $u>u_\mu$ in $\Sigma_\mu$. On the other 
hand, we can also exclude that $\bar x\in T_\mu$; indeed, we have 
$$
 u(x_n)-u((x_n)_{\la_n})= 2(x_n- \la_n)\myfrac{\partial u}{\partial 
x_1} (\xi_n)
$$
for a point $\xi_n\in ((x_n)_{\la_n}, x_n)$. If (a  subsequence of) 
$x_n$ converges to  a point in $T_\mu$, then for $n$ large we have 
dist$(\xi_n, T_\mu)<\sigma$ and from \rife{sig} we get 
$u(x_n)-u((x_n)_{\la_n})>0$ contradicting \rife{cont}.   We are left 
with the possibility that $\bar x\in \partial \Sigma_\mu \setminus 
T_\mu$: but this is also a contradiction since  $u$ blows up at the 
boundary and  it is locally bounded in the interior, so that  
$u(x_n)-u((x_n)_{\la_n})$ would converge to infinity.

Thus $\mu=0$ and \rife{lam} holds in the whole $\{x\in B_R\,:\, 
x_1>0\}$. We deduce  that $u$ is symmetric in the $x_1$ direction
and $\partial u/\partial x_1 >0$. Applying to any other direction we 
conclude that $u$ is radial and 
$\partial u/\partial r>0$.    
\qed

\begin{remark}{\rm
Let us recall  that in some special  examples (for instance when $g(s)$ has an exponential or  a power--like growth) the asymptotic behavior at the boundary of the gradient of the large solutions has already been studied (see e.g.  \cite{BE}, \cite{BM1}, \cite{PV}) so that the previous result could be directly applied to prove symmetry. In general, through a blow--up argument, we are able to prove \rife{der} if 
\begin {equation}\label{b-u}
s\mapsto \frac{g(s)}{\sqrt{G(s)}}\int_s^{\infty}\frac1{\sqrt{2G(\xi)}}d\xi
\end {equation}
 is bounded at infinity; however this assumption does not include the case when $g$ has a  slow growth at infinity (such as $g(r)\equiv r(\ln r)^\ga$ with $\ga>2$) and
is not so general as  \rife{g2}.}
\end{remark}

\vskip1em

\begin{theorem}\label{main}
Assume that \rife{g1}, \rife{g2} and \rife{g3} hold. 
Then any solution $u$ of \rife{large2} is radial and $\prt u/\prt 
r>0$ in $B_{R}\setminus \{0\}$.
\end{theorem}

The following preliminary result is  a consequence of more general 
results in \cite{MV}, \cite{MV1}. However we provide here a simple 
self-contained proof  for the radial case.

\begin{lemma}\label{unico}
Let $h$ be a  convex increasing function satisfying the 
Keller-Osserman condition
\begin{equation}\label{ok}
\int_{a}^{+\infty} \myfrac{ds}{\sqrt {H(s)}}<\infty\,,\qquad 
H(s)=\int_a^s h(t)dt\,,
\end{equation}
for some $a>0$. Then the problem
\begin{equation}\label{h}
\cases{
-\Delta v+h(v)=0 & in $B_R$\cr
\noalign{\medskip} 
\;\;\;\,\lim\limits_{|x|\to R}v(x)=\infty ,& \cr}
\end{equation}
has a unique solution.
\end{lemma}
\proof  Since $h$ is increasing, there exist a maximal and a minimal 
solution $\overline v$ and $\underline v$, which are both radial, so 
that it is enough to prove that $\overline v=\underline v$. To this 
purpose, observe that if $v$ is radial we have
$(v' r^{N-1})'=r^{N-1}h(v)$ so that, since $v'(0)=0$, and 
replacing $H$ by $\tilde H=H-H(\min v)$ which is nonnegative on the 
range of values of $v$, we have
$$
\myfrac{(v' r^{N-1})^2}2=\int_0^r  s^{2(N-1)}h(v)v' ds \leq  
r^{2(N-1)}\tilde H(v)
$$ 
which yields 
\begin{equation}\label{co}
0\leq v'< \sqrt{2\tilde H(v)}\,.
\end{equation}
Define now $w=F(v)=\myint{v}{\infty}\myfrac{ds}{\sqrt{2\tilde 
H(s)}}$.  A straightforward computation, and  condition \rife{ok},   
show that $w$ solves the problem
\begin{equation}\label{ew}
\cases{\Delta w= b(w)\left( |Dw|^2-1\right) \,& in $B_R$,\cr 
\;\;\;w=0& on $\partial B_R$,\cr}
\end{equation}
where $b(w)=h(v)/\sqrt{2\tilde H(v)}$. One can easily check that the 
convexity of $h$ implies that 
$h(v)/\sqrt{2\tilde H(v)}$ is nondecreasing, hence $b(w)$ is 
nonincreasing with respect to $w$. Moreover, since 
$|Dw|=|w'|=v'/\sqrt{2\tilde H(v)}$, from \rife{co}  one gets $|w'|<1$. 
Note that the transformation $v\mapsto w$ establishes a one-to-one monotone 
correspondence between the large solutions of \rife {h} and 
the solutions of \rife {ew}, so that $\overline w=F(\underline v)$ 
and $\underline w=F(\overline v)$  are respectively the  minimal and 
the maximal solutions
of  \rife{ew}. Thus we have 
$$
((\overline w-\underline w)'r^{N-1})' = r^{N-1}\left[  b(\overline w) 
(|\overline w' |^2-1) - b(\underline w) (|\underline w' |^2-1)\right] 
\geq  r^{N-1} b(\overline w) (|\overline w' |^2-|\underline w'|^2)\,, 
$$
so that the function $z= (\overline w-\underline w)'r^{N-1}$ satisfies
$$
z' \geq a(r) z\,,\quad z(0)=0\,,
\qquad \hbox{where $a(r)=  b(\overline w) (\overline w' +\underline 
w')$.}
$$
Because $a$ is locally bounded on $[0,R)$, we deduce that $z\geq 0$, 
hence $\overline w-\underline w$ is nondecreasing. Since $\overline 
w-\underline w$ is nonnegative and $\overline w(R)=\underline w(R)=0$ 
we deduce that $\overline w=\underline w$, hence $\overline 
v=\underline v$.
\qed

\vskip2em

\begin{lemma}\label{d}
Assume that $g$ satisfies (\ref{g2}) and \rife{g3}, and that $u$ is a 
solution of (\ref {large2}). Then
\begin{equation}\label{infde}
\BA{l} 
(i)\quad\displaystyle\lim_{|x|\to R}\nabla_{\tau}u(x)= 0\\ 
[3mm]
(ii)\quad\displaystyle\lim_{|x|\to R}\myfrac{\partial u}{\partial r}(x)=\infty,
\EA
\end{equation}
and the two limits hold uniformly with respect to $\{x:\abs x=r\}$.
\end{lemma}
\proof  In spherical coordinates $(r,\gs)\in (0,\infty)\ti S^{N-1}$ the Laplace operator 
takes the form
$$
\Delta u = \myfrac{\partial^2 u}{\partial r^2} + \myfrac {N-1}r 
\myfrac{\partial u}{\partial r}+ \myfrac 1{r^2} \Delta_s u,
$$
 where $\Delta_s $ is the Laplace Beltrami operator on $S^{N-1}$. If 
 $\{\gg_{j}\}_{j=1}^{N-1}$ is a system of $N-1$ geodesics on 
 $S^{N-1}$ crossing orthogonally   
 at $\tilde\gs$, there holds
\begin {equation}\label {LB}
\Delta_s u(r,\tilde\gs)=\sum_{j\geq 1}{}\myfrac 
 {d^{2}u(r,\gg_{j}(t))}{dt^{2}}|_{t=0}.
\end{equation}
 On the sphere the geodesics are large circles. The system of 
 geodesics can be obtained by considering a set of skew symmetric 
 matrices $\{A_{j}\}_{j=1}^{N-1}$ such that $\langle A_{j}\tilde\gs,
 A_{k}\tilde\gs\rangle =\gd_{j}^k$, and by putting $\gg_{j}(t)=e^{tA_{j}}\tilde\gs$.
\smallskip

\nind{\it Step 1: two-side estimate on the tangential first derivatives}. By 
assumption \rife{g2} $g$ can be written as
$$
g(s)= g_\infty(s) + \tilde g(s),
$$
where $g_\infty(s)$ is a  convex increasing function satisfying 
\rife{g3} and $\tilde g(s)$ is  a locally Lipschitz function such 
that $\tilde g\equiv 0$ in $[M,\infty)$ for some $M>0$.
In particular, $u$ satisfies
$$
\Delta u - g_\infty(u) = \tilde g(u).
$$
Since $u$ blows up uniformly, there holds $u(x)\geq M$ if $\abs 
x\in [r_0,R)$ for a  certain $r_0<R$, hence 
$$
\tilde g(u)=\tilde g(u)\chi_{\{\abs x\leq r_0\}}\leq K_0, 
$$
so that 
\begin{equation}\label{lap}
|\Delta u - g_\infty(u) |\leq K_0\,.
\end{equation}
Set $\vfi(r)=\myfrac1{2N}(R^2-r^2)$, thus $\vfi$ satisfies $-\Delta 
\vfi=1$ and $\vfi=0$ on $\partial B_R$. 
We deduce from \rife{lap} that 
$$
\Delta (u-K_0 \vfi)\geq g_\infty(u)\geq g_\infty(u-K_0\vfi),
$$
since $g_\infty$ is increasing and $\vfi$ is nonnegative. Thus 
$u-K_0\vfi$ is a sub-solution of the problem 
\begin{equation}\label{infi}
\cases{
-\Delta v+g_\infty(v)=0 & in $B_R$\cr
\noalign{\medskip} 
\qquad\!\!\lim\limits_{|x|\to R}v(x)=\infty.& \cr}
\end{equation}
Similarly $u+K_0\vfi$ is a super-solution of the same problem. By 
Lemma \ref{unico}, problem \rife{infi} has a unique solution $U_R$. 
By approximating $U_{R}$ by the large solution $U_{R'}$ of the same 
equation in $B_{R'}$ with $R'<R$ and $R'>R$ we derive
\begin{equation}\label{U}
U_R- K_0\vfi\leq u\leq U_R+K_0 \vfi\,.
\end{equation}
\vskip1em

Since the problem (\ref {large2}) is invariant by rotation, for any 
$j=1,\ldots,N-1$, and any $h\in\BBR$, the function $u^h$ defined by
$u^h(x)=u(e^{hA_{j}}(x))=u(r,e^{hA_{j}}\gs)$ is a solution of (\ref {large2}) and still 
$\tilde g(u^h)=0$ if $r\in [r_0,R)$, so that
$$
\Delta u^h= g_\infty(u^h) \qquad \hbox{if $r\in (r_0,R)$.}
$$
Since  $u\in C^1(B_{R})$, there holds
$$
|u^h-u| \leq L \,|h|  \qquad \hbox{if $r=r_0$.}
$$
Let us set 
\begin{equation}\label{pot}
P(r)=\cases{\myfrac{r^{2-N}-R^{2-N}}{r_{0}^{2-N}-R^{2-N}}&  if 
$N>2$\cr \noalign{\medskip}
 \myfrac{\ln r-\ln R}{\ln r_{0}-\ln R}&  if $ N=2$,\cr}
\end{equation}
and $v^h(x)=u^h(x)+|h| L P(\abs x)$;  then $\Delta v^h=\Delta u^h$, 
and since $g_\infty$ is increasing,
\begin{equation}\label{xx}
\Delta (v^h- u) \leq g_\infty(v^h)- g_\infty(u)\qquad\mbox{in 
}B_{R}\setminus B_{r_{0}}.
\end{equation}
Observe that  $u^h$, as $u$, also satisfies \rife{U}, so that in particular 
\begin{equation}\label{tb}
u^h(x) - u(x)\to 0\quad\hbox{as $\abs x\to R$.}
\end{equation}
Therefore $v^h(x)-u(x)\to 0$  as $\abs x\to 
R$ too, while by construction $v^h\geq u$ on $\partial B_{r_0}$. We conclude from \rife{xx} (e.g.  using the test function $(v^h-u+\vep)_{-}$, which is compactly supported, and then letting $\vep$ go to zero) that 
$$
v^h= u^h + |h| L P(r) \geq u. 
$$
We recall that the Lie derivative $L_{A_{j}}u$ of $u(r,.)$  following the 
vector field tangent to $S^{N-1}$ $\eta\mapsto A_{j}\eta$  is defined by
$$L_{A_{j}}u(r,\gs)=\myfrac {du(r,e^{tA_{j}}\gs)}{dt}|_{t=0},
$$
so we get, by letting $h\to 0$,
\begin{equation}\label{pri}
\left|L_{A_{j}}u(r, \tilde \gs)\right| \leq  L P(r) < C (R-r)\,.
\end{equation}
\smallskip

\nind{\it Step 2: one-side estimate on the tangential second derivatives}.
Next we define  the function $w^h$ by
$$
w^h= \myfrac{u^h+u^{-h}-2u}{h^2}\,.
$$
As before, let $r_0<R$ be such that $u\geq M$ on $B_{R}\setminus 
B_{r_{0}}$. Thus $g(u)=g_\infty(u)$ on $B_{R}\setminus B_{r_{0}}$,
and
$$
\Delta w^h= \myfrac1{h^2} \left(g_\infty(u^h)+ g_\infty(u^{-h})-2 
g_\infty(u)\right)\,
\qquad \mbox{on }B_{R}\setminus B_{r_{0}}.
$$
Since $g_\infty$ is convex, there holds
$$
g_\infty(a)+g_\infty(b)-2g_\infty(c)\geq \xi (a+b-2c)\qquad\forall 
\xi\in \partial g_\infty(c)
$$
where $\partial g_\infty(c)=[g'_\infty(c_{-}),g'_\infty(c_{+})]$. 
Hence
$$
\Delta w^h \geq \xi_u\, w^h\qquad \mbox{  on }B_{R}\setminus 
B_{r_{0}},
$$
for any $\xi_u\in \partial g_\infty (u)$. Since  $g_\infty$ is increasing, we have $\xi_u\geq 0$, 
therefore $(w^h)_+$ is a 
subharmonic function in $B_R\setminus B_{r_0}$. 
As $u\in C^2(B_R)$, 
there exists  $\tilde L>0$ such that 
$$
w^h\leq \tilde L\qquad \hbox{on $\partial B_{r_0}$.}
$$
Moreover from \rife{pri} we get that $w^h=0$ on $\partial B_R$. We 
conclude that 
$$
(w^h)_+(x) \leq  \tilde L \,P(\abs x) \,,
$$
where $P(r)$ is defined in \rife{pot}. Letting $h$ tend to zero we 
obtain
\begin{equation}\label{sec}
\myfrac {d^{2}u(r,e^{tA_{j}}\tilde\gs)}{dt^{2}}|_{t=0} \leq  \tilde L 
\,P(r)\qquad \mbox{  for }r\in [r_{0},R ).
\end{equation}
Using (\ref {LB}), and the fact that $\tilde\gs$  is arbitrary, we 
derive
\begin{equation}\label{sec'}
\Delta_s u(r, \gs)\leq (N-1)\tilde L \,P(r)\forevery (r,\gs)\in 
[r_{0},R)\ti S^{N-1}.
\end{equation}
 
\nind{\it Step 3: estimate on the radial derivative}. Using \rife {pot} 
and \rife{sec'} we 
deduce that
$$
(\Delta_s u )^+ (x)=\circ (1)\quad \hbox{uniformly as $\abs x\to R$.}
$$ 
Therefore
\begin{equation}\label{ra}
\begin{array}{rl}
\myfrac{\partial}{\partial r}\left(r^{N-1}\myfrac{\prt  u}{\prt 
r}\right) &= r^{N-1}\left[g(u)-\myfrac1{r^2}\Delta_s u\right]
\\ 
\noalign{\medskip}
&
\geq r^{N-1}g_\infty(u)-\circ(1)\quad \hbox{uniformly as $r\to R$.}
\end{array}
\end{equation}
Now one can easily conclude: let $z(r)$ denote the minimal (hence 
radial) solution of 
$$
\left\{\BA {l}
\Delta z= g_\infty(z) \quad \mbox { in }B_R\setminus B_{r_0}\\[2mm] 
\phantom {\Delta}
z=\min\limits_{\partial B_{r_0}} u  \quad \mbox { on }\partial 
B_{r_0}\\\phantom {\Delta}
 \lim\limits_{r\to R}z=\infty.
 \EA\right.
$$
We have $u(x)\geq z(x)$ if $\abs x\in [r_0,R)$, hence 
$g_\infty(u)\geq g_\infty(z)$. Because this last function is not 
integrable near $\prt B_{R}$, one obtains
$$
\lim\limits_{r\to R} \quad \int_{r_0}^r s^{N-1}g_\infty(u(s,\gs))ds\to 
\infty \quad \hbox{uniformly for $\gs\in S^{N-1}$.}
$$
Clearly \rife{ra} implies
$$
 \pd ur(r,\gs)\,\,\mathop{\to}^{r\to R}\, \infty\quad \hbox{uniformly for $\gs\in S^{N-1}$}.
$$
This completes the proof of \rife{infde}.\qed
\medskip
\vskip1em
\nind {\bf Proof of Theorem \ref{main}.}
By assumptions \rife{g2} and \rife{g3}, and Lemma \ref{d}, we  deduce 
that $u$ satisfies \rife{der}, hence we apply Lemma \ref{GNN} to 
conclude.
\qed
\vskip2em
Finally,  let us point out that thanks to Lemma \ref{d} and using the moving plane method as in  Theorem \ref{GNN}, 
we can derive a result describing the boundary behaviour of any solution of
\begin {equation}\label {E1}\left\{\BA {l}
-\Gd u+g(u)=0\quad \mbox {in } \Gg_{R,r}=\{x\in\BBR^N:r<\abs x<R\}\\\phantom{-,}
\displaystyle{\lim_{\abs x\to R}}u(x)=\infty, 
\EA\right.
\end {equation}
which extends  a similar result in \cite{GNN}.

\bcor {torus} Assume that $g$ satisfies \rife{g1}, \rife{g2} and \rife{g3}. Then any solution of  (\ref{E1}) satisfies \rife{infde}
and verifies $\prt_{r}u>0$ on $\Gg_{R,(R+r)/2}$.
\es

\end{document}